\documentclass {article}

\usepackage {amsmath}
\usepackage {amssymb}
\usepackage {color}
\usepackage {graphicx}
\usepackage {oldgerm}

\newcommand {\Aut}  {\ensuremath{\operatorname{Aut}}}
\newcommand {\cn}   {\ensuremath{\operatorname{c}}}
\newcommand {\cp}   {\ensuremath{\mathbin{\square}}}
\newcommand {\D}    {\ensuremath{\operatorname{D}}}
\newcommand {\id}   {\ensuremath{\text{\rm id}}}
\newcommand {\m}    {\ensuremath{\operatorname{m}}}
\newcommand {\Prob} {\ensuremath{\operatorname{Pr}}}
\newcommand {\supp} {\ensuremath{\operatorname{supp}}}
\newcommand {\V}    {\ensuremath{\operatorname{V}}}

\newcommand {\prf}  {\noindent \textit{Proof.} }
\newcommand {\qed}  {\hfill $\Box$}
\newcommand {\qedd} {\hfill $\Box$\medskip}

\newtheorem {lem}       {Lemma}     [section]
\newtheorem {thm} [lem] {Theorem}
\newtheorem {cor} [lem] {Corollary}

\numberwithin {equation} {section}

\title  {Distinguishing graphs with \\ infinite motion and nonlinear growth}
\author {Johannes Cuno\thanks{All three authors, namely Johannes Cuno, Wilfried Imrich, and Florian Lehner, acknowledge the support of the Austrian Science Fund (FWF), project W1230-N13.} \\ Technische Universit\"{a}t Graz \\ \texttt{cuno@math.tugraz.at} \and Wilfried Imrich \\ Montanuniversit\"at Leoben \\ \texttt{imrich@unileoben.ac.at} \and Florian Lehner \\ Technische Universit\"{a}t Graz \\ \texttt{f.lehner@tugraz.at}}

\begin{document}

\maketitle

\begin{abstract}
The distinguishing number $\D(G)$ of a graph $G$ is the least cardinal $d$ such that $G$ has a labeling with $d$ labels which is only preserved by the trivial automorphism. We show that the distinguishing number of infinite, locally finite, connected graphs $G$ with infinite motion and growth $o(n^2/\log_2 n)$ is either $1$ or $2$, which proves the Infinite Motion Conjecture of Tom Tucker for this type of graphs. The same holds true for graphs with countably many ends that do not grow too fast. We also show that graphs $G$ of arbitrary cardinality are $2$-distinguishable if every nontrivial automorphism moves at least uncountably many vertices $\m(G)$, where ${\m(G)} \geq |\!\Aut(G)|$. This extends a result of Imrich et al.~to graphs with automorphism groups of arbitrary cardinality. \smallskip

% ================================================================================
% STICHWORTE
% ================================================================================

\noindent \textbf{Keywords:} Distinguishing number, automorphisms, infinite graphs. \smallskip

% ================================================================================
% AMS-KLASSIFIKATION
%
% 05C25 - Graphs and abstract algebra (groups, rings, fields, etc.)
% 05C63 - Infinite graphs
% 05C15 - Coloring of graphs and hypergraphs
% 03E10 - Ordinal and cardinal numbers
% ================================================================================

\noindent \textbf{MSC classes:} 05C25, 05C63, 05C15, 03E10.
\end{abstract}

\section{Introduction}
Albertson and Collins \cite{alco-96} introduced the distinguishing number $\D(G)$ of a graph $G$ as the least cardinal $d$ such that $G$ has a labeling with $d$ labels which is only preserved by the trivial automorphism.

This seminal concept spawned many papers on finite and infinite graphs. We are mainly interested in infinite, locally finite, connected graphs of polynomial growth, see \cite{imjekl-08}, \cite{tu-06}, \cite{smtuwa-xx}, and in graphs of higher cardinality, see \cite{imkltr-07}, \cite{langu-10}. In particular, there is one conjecture on which we focus our attention, the Infinite Motion Conjecture of Tom Tucker.

Before stating it, we introduce the notation $\m(\phi)$ for the number of elements moved by an automorphism $\phi$, and call $\m(\phi)$ the \emph{motion} of $\phi$. In other words, $\m(\phi)$ is the size of the set of vertices which are not fixed by $\phi$, that is, the size of its \emph{support}, $\supp(\phi)$.

\medskip
\textbf{The Infinite Motion Conjecture of Tom Tucker.} \emph{Let $G$ be an infinite, locally finite, connected graph. If every nontrivial automorphism of $G$ has infinite motion, then the distinguishing number $\D(G)$ of $G$ is either 1 or 2.}
\medskip

For the origin of the conjecture and partial results compare \cite{smtuwa-xx}. The conjecture is true if $\Aut(G)$ is countable, hence we concentrate on graphs with uncountable group.

The validity of the conjecture for graphs with countable group follows from either one of two different results in \cite{imsmtuwa-xx}. One of them replaces the requirement of infinite motion by a lower and upper bound on the size of the automorphism group. It asserts that every infinite, locally finite, connected graph $G$ whose automorphism group is infinite, but strictly smaller than $2^{\aleph_0}$, has countable group, infinite motion, and distinguishing number 2. For a precise formulation see Theorem \ref{thm:simon}. The proof is not easy and follows from results of either Halin \cite{ha-73}, Trofimov \cite{tr-85}, or Evans \cite{ev-87}.

The other one relaxes the condition of local finiteness and requires that the group is at most countable. It asserts that countably infinite, connected graphs with finite or countably infinite group and infinite motion are $2$-distinguishable, no matter whether they are locally finite or not, see Theorem \ref{thm:tucker}. The proof is short and elementary.

For uncountable connected graphs with countable motion the Infinite Motion Conjecture need not be true. We turn to this case in Section \ref{sec:generalization}, suggest a version of the conjecture for uncountable connected graphs, and prove its validity under a bound on the size of the automorphism group.

\section{Preliminaries}
Throughout this paper the symbol $\mathbb{N}$ denotes the set $\{1,2,3,\ldots\}$ of positive integers, whereas the symbol $\mathbb{N}_{0}$ refers to the set $\{0,1,2,3,\ldots\}$ of non-negative integers.

Let $G$ be a graph with vertex set $\V(G)$. Let $X$ be a set. An $X$-\emph{labeling} $l$ of $G$ is a mapping $l:\V(G)\rightarrow X$. For us $X$ will mostly be the set $\{\text{black},\text{white}\}$. In this case, we speak of a \emph{$2$-coloring} of $G$.

Let $l$ be an $X$-labeling of $G$. Consider an automorphism $\phi\in\Aut(G)$. If, for every $v\in\V(G)$, $l({\phi(v)}) = l(v)$, we say that $l$ \emph{is preserved by} $\phi$. {If this is not the case}, we say that $l$ \emph{breaks} $\phi$. An $X$-labeling $l$ of $G$ is called \emph{distinguishing} if it is only preserved by the trivial automorphism. The \emph{distinguishing number} $\D(G)$ of $G$ is the least cardinal $d$ such that there exists a distinguishing $X$-labeling of $G$ with $|X|=d$.

Given a group $A$ equipped with a homomorphism $\phi:A\rightarrow\Aut(G)$, we say that $A$ \emph{acts on} $G$. Moreover, we say that $A$ \emph{acts nontrivially} on $G$ if there is an $a \in A$ such that $\phi(a)$ moves at least one vertex of $G$. By abuse of language we write $a(v)$ instead of $\phi(a)(v)$ and say that an $X$-labeling $l$ of $G$ \emph{is preserved by} $a\in A$ if it {is preserved by} $\phi(a)\in\Aut(G)$.

The \emph{ball with center $v_0\in\V(G)$ and radius $r$} is the set of all vertices $v\in\V(G)$ with $d_G(v_0,v)\leq r$ and is denoted by $B_{v_0}^{G}(r)$, whereas $S_{v_0}^{G}(r)$ stands for the set of all vertices $v\in\V(G)$ with $d_G(v_0,v)=r$. We call it the \emph{sphere with center $v_0\in\V(G)$ and radius $r$.} If $G$ is clear from the context, we just write $B_{v_0}(r)$ and $S_{v_0}(r)$ respectively. For terms not defined here we refer to \cite{haimkl-11}.

Although our graphs are infinite, as long as they are locally finite, all balls and spheres of finite radius are finite. The number of vertices in $B_{v_0}^{G}(r)$ is a monotonically increasing function of $r$, because
\[\left|B_{v_0}^{G}(r)\right| = \sum_{i={0}}^r\left|S_{v_0}^{G}(i)\right|\quad\text{and}\quad\left|S_{v_0}^{G}(i)\right| \geq 1\,.\]
Nonetheless, the growth of $\left|B_{v_0}^{G}(r)\right|$ depends very much on $G$, and it is helpful to define the growth rate of a graph.

We say that an infinite, locally finite, connected graph $G$ has \emph{polynomial growth} if there is a vertex $v_0 \in \V(G)$ and a polynomial $p$ such that
\[\forall\,r\in\mathbb{N}_{0}\colon\left|B_{v_0}^{G}(r)\right| \leq p(r)\,.\]
It is easy to see that this implies that all functions $\left|B_{v}^{G}(r)\right|$ are bounded by polynomials of the same degree as $p$, independent of the choice of $v \in \V(G)$. In this context it should be clear what we mean by \emph{linear} and \emph{quadratic growth}. Observe that the two-sided infinite path has linear growth, and that the growth of the grid of integers in the plane is quadratic.

We say that $G$ has \emph{exponential growth} if there is a constant $c>1$ such that
\[\forall\,r\in\mathbb{N}_{0}\colon\left|B_{v_0}^{G}(r)\right| \geq c^r\,.\]
Notice that homogeneous trees of degree $d > 2$, that is, infinite trees where every vertex has the same degree $d$, have exponential growth. For the distinguishability of such trees and tree-like graphs, see \cite{wazh-07} and \cite{imkltr-07}.

We are mainly interested is the distinguishability of infinite, locally finite, connected graphs of polynomial growth. For us, the following lemma will be helpful.

\begin{lem}
\label{lem:halfzero}
Let $A$ be a finite group acting on a graph $G$. If a coloring of $G$ breaks some element of $A$, then it breaks at least half of the elements of $A$.
\end{lem}

\prf The elements of $A$ that preserve a given coloring form a subgroup. If some element of $A$ is broken, then this subgroup is proper and thus, by Lagrange's theorem, cannot contain more than half of the elements of $A$. \qedd

If the action is nontrivial, then we can always find a coloring that breaks at least one element. Hence, we have the following result.

\begin{lem}
\label{lem:halfone}
Let $G$ be a graph. If $A$ is a finite group acting nontrivially on $G$, then there exists a $2$-coloring of $G$ that breaks at least half of the elements of $A$.
\end{lem}

The proof of Lemma \ref{lem:halfone} is based on the fact that $A$ is a group. But a very similar result holds for any finite family of nontrivial automorphisms, as the following lemma shows.

\begin{lem}
\label{lem:halftwo}
Let $G$ be a finite graph. If $A$ is a finite set equipped with a mapping $\phi:A\rightarrow\Aut(G)\setminus\{\id\}$, then there exists a $2$-coloring of\, $G$ that breaks $\phi(a)$ for at least half of the elements of $A$.
\end{lem}

\prf Let $\V(G)=\{v_{1},v_{2},\ldots,v_{n}\}$. For every $k\in\{1,2,\ldots,n\}$, let $A_{k}$ be the set of all $a\in A$ with $\supp(\phi(a)) \subseteq\{v_{1},v_{2}, \ldots,v_{k}\}$. We show by induction that the assertion holds for all $A_{k}$ and, in particular, for $A$. Because $A_{1}$ is the empty set, the assertion is true for $A_{1}$. Suppose it is true for $A_{k-1}$. Then we can choose a $2$-coloring of $G$ that breaks $\phi(a)$ for at least half of the elements of $A_{k-1}$. This remains true, even when we change the color of $v_{k}$. Notice that, for every $a\in A_{k}\setminus A_{k-1}$, $\phi(a)$ either maps $v_{k}$ into a white vertex in $\{v_{1},v_{2},\ldots,v_{k-1}\}$ or into a black vertex in $\{v_{1},v_{2},\ldots,v_{k-1}\}$. Depending on which of the two alternatives occurs more often, we color $v_{k}$ black or white such that this $2$-coloring also breaks $\phi(a)$ for at least half of the elements of $A_{k}\setminus A_{k-1}$ and, hence, for at least half of the elements of $A_{k}$. \qedd

If every nontrivial automorphism of a graph $G$ has infinite motion, we say that $G$ has \emph{infinite motion}. For such graphs the following result from \cite{imsmtuwa-xx} will be of importance.

\begin{lem}
\label{lem:main}
Let $G$ be an infinite, locally finite, connected graph with infinite motion. If an automorphism $\phi\in\Aut(G)$ fixes a vertex $v_0\in\V(G)$ and moves at least one vertex in $S_{v_0}(k)$, then, for every $i\geq k$, it moves at least one vertex in $S_{v_0}(i)$.
\end{lem}

\section{Graphs of nonlinear growth}
\label{sec:nonlinear}
In \cite{imsmtuwa-xx}, it was shown that infinite, locally finite, connected graphs with infinite motion and linear growth have countable automorphism group, and therefore distinguishing number either 1 or 2.

If the growth rate of such graphs becomes nonlinear, then the automorphism group can become uncountable. This holds, even if the growth rate becomes only slightly nonlinear.

\begin{thm}
\label{thm:stretchedtree}
Let $\varepsilon>0$. Then there exists an infinite, locally finite, connected graph $G$ with uncountable automorphism group, infinite motion, and nonlinear growth function $g:\mathbb{N}_{0}\rightarrow\mathbb{N}_{0}$ such that, for sufficiently large $n\in\mathbb{N}_{0}$, $g(n)$ is bounded from above by $n^{1+\varepsilon}$.
\end{thm}

\prf We construct $G$ from $T_3$, that is, the tree in which every vertex has degree 3. First, choose an arbitrary vertex $v_0\in\V(T_3)$. Our strategy is to replace the edges of $T_3$ by paths such that, for sufficiently large $n\in\mathbb{N}_{0}$, $g(n)=\left|B_{v_0}^G(n)\right|\leq n^{1+\varepsilon}$.

For every $i\in\mathbb{N}_{0}$, there are $3\cdot 2^{i}$ edges from $S_{v_0}^{T_3}(i)$ to $S_{v_0}^{T_3}(i+1)$. If we replace them by paths of the same length, then the cardinality of the balls $B_{v_0}^{G}(n)$ grows linearly with slope $3\cdot 2^{i}$ from $S_{v_0}^{T_3}(i)$ to $S_{v_0}^{T_3}(i+1)$.

Observe that, given any affine linear function $h:\mathbb{N}_{0}\rightarrow\mathbb{N}_{0}$, there is a number $n_h\in\mathbb{N}$ such that, for all $n \geq n_h$, $h(n)\leq n^{1+\varepsilon}$. In particular, we may consider the functions $h_i:\mathbb{N}_{0}\rightarrow\mathbb{N}_{0}$ defined by $h_i(x) = 3\cdot 2^i\cdot x+1$, and choose numbers $n_i\in\mathbb{N}$ such that, for every $n \geq n_i$, $h_i(n)\leq n^{1+\varepsilon}$.

\begin{figure}
\begin{center}
\scalebox{1}{\input{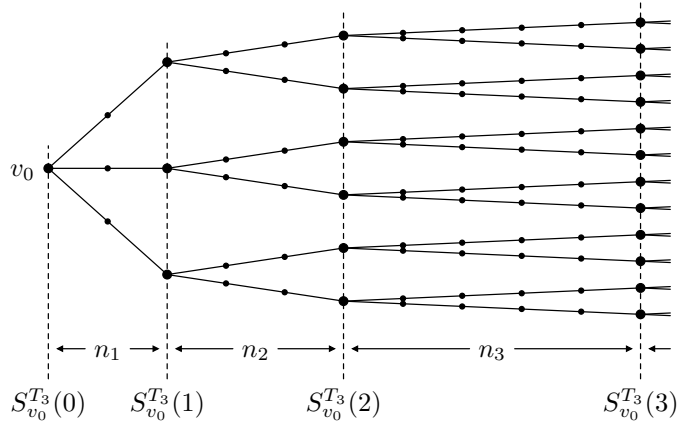}} \vspace{-2mm}
\end{center}
\caption{Replacing the egdes of $T_3$ by paths.}
\label{fig:replacing}
\end{figure}

As illustrated in Figure \ref{fig:replacing}, for every $i\in\mathbb{N}_{0}$, we replace the edges from $S_{v_0}^{T_3}(i)$ to $S_{v_0}^{T_3}(i+1)$ by paths of length $n_{i+1}$. For every $i\in\mathbb{N}$ and every vertex $v\in\V(G)$ on such a path from $S_{v_0}^{T_3}(i)$ to $S_{v_0}^{T_3}(i+1)$, we have $d_G(v,v_0)\geq n_{i}$ and, hence,
\[ g(d_G(v,v_0))\leq 3\cdot 2^i\cdot d_G(v,v_0)+1=h_{i}(d_G(v,v_0))\leq d_G(v,v_0)^{1+\varepsilon}\,. \]
So, for every $n \geq n_1$, $g(n)$ is bounded from above by $n^{1+\varepsilon}$. Every automorphism of $T_3$ that fixes $v_0$ induces an automorphism of $G$. It is easy to see that this correspondence is bijective. Thus, $\Aut(G)$ is uncountable. Furthermore, $G$ inherits infinite motion from $T_3$. Since $\Aut(G)$ is uncountable, the result of \cite{imsmtuwa-xx} mentioned at the beginning of Section \ref{sec:nonlinear} implies that $G$ cannot have linear growth. \qedd

Though we cannot assume that the automorphism groups of our graphs are countable, we prove that infinite, locally finite, connected graphs with infinite motion and nonlinear, but moderate, growth are still $2$-distinguishable, that is, they have distinguishing number either $1$ or $2$.

Our construction of a suitable coloring consists of several steps. In Lemma \ref{lem:fixroot} we color a part of the vertices in order to break all automorphisms that move a distinguished vertex $v_{0}$. In Lemma \ref{lem:breakautos} we show how to color some of the remaining vertices in order to break more automorphisms. Iteration of this procedure yields a distinguishing coloring, as shown in Theorem \ref{thm:spheres}.

\begin{lem}
\label{lem:fixroot}
Let $G$ be an infinite, locally finite, connected graph with infinite motion and $v_0\in\V(G)$. Then, for every $k\in\mathbb{N}$, one can $2$-color all vertices in $B_{v_0}(k+3)$ and $S_{v_0}(\lambda k+4)$, $\lambda\in\mathbb{N}$, such that, no matter how one colors the remaining vertices, all automorphisms that move $v_{0}$ are broken.
\end{lem}

\begin{figure}
\begin{center}
\scalebox{1}{\input{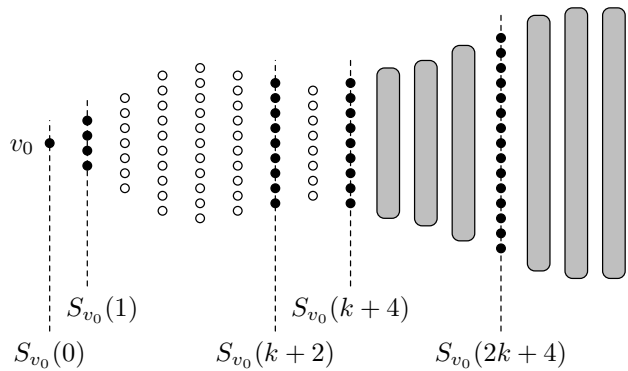}} \vspace{-2mm}
\end{center}
\caption{Breaking all automorphisms that move $v_{0}$.}
\label{fig:fixroot}
\end{figure}

\prf If $k=1$, then we color $v_0$ black and all $v\in\V(G)\setminus\{v_0\}$ white, whence all automorphisms that move $v_{0}$ are broken. So, let $k\geq 2$. First, we color all vertices in $S_{v_0}(0)$, $S_{v_0}(1)$, and $S_{v_0}(k+2)$ black and the remaining vertices in $B_{v_0}(k+3)$ white. Moreover, we color all vertices in $S_{v_0}(\lambda k+4)$, $\lambda\in\mathbb{N}$, black and claim that, no matter how we color the remaining vertices, $v_{0}$ is the only black vertex that has only black neighbors and only white vertices at distance $r\in\{2,3,\ldots, k+1\}$, see Figure \ref{fig:fixroot}.

It clearly follows from this claim that this coloring breaks every automorphism that moves $v_0$. It only remains to verify the claim.

Consider a vertex $v\in\V(G)\setminus\{v_0\}$. If $v$ is not in $S_{v_0}(1)$, then it is easy to see that $v$ cannot have the aforementioned properties. So, let $v$ be in $S_{v_0}(1)$ and assume it has only black neighbors and only white vertices at distance $2$. Then it cannot be neighbor to any vertex in $S_{v_0}(2)$, but must be neighbor to all vertices in $B_{v_0}(1)$ except itself. Therefore, the transposition of the vertices $v$ and $v_0$ is a nontrivial automorphism of $G$ with finite support. Since $G$ has infinite motion, this is not possible. \qed

\begin{lem}
\label{lem:breakautos}
Let $G$ be an infinite, locally finite, connected graph with infinite motion and $v_0\in\V(G)$. Moreover, let $\varepsilon>0$. Then there exists a $k\in\mathbb{N}$ such that, for every $m\in\mathbb{N}$ and for every $n\in\mathbb{N}$ that is sufficiently large and fulfills
\begin{equation}
\label{eq:1}
\left|S_{v_0}(n)\right| \,\leq\, \frac{n}{(1+\varepsilon)\log_2 n }\,,
\end{equation}
one can $2$-color all vertices in $S_{v_0}(m+1),S_{v_0}(m+2),\ldots,S_{v_0}(n)$, but not those in $S_{v_0}(\lambda k+4)$, $\lambda\in\mathbb{N}$, such that all automorphisms that fix $v_{0}$ and act nontrivially on $B_{v_0}(m)$ are broken.%\footnote{See Figure \ref{fig:breakautos} for illustration.}
\end{lem}

\begin{figure}
\begin{center}
\scalebox{1}{\input{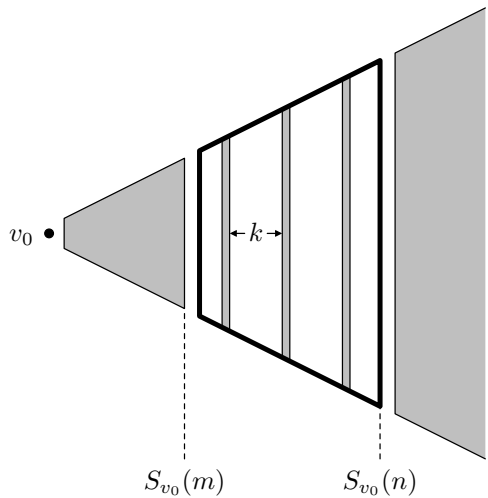}} \vspace{-2mm}
\end{center}
\caption{Breaking all automorphisms that fix $v_{0}$ and act nontrivially on $B_{v_0}(m)$.}
\label{fig:breakautos}
\end{figure}

The coloring and the meaning of the variables $m$, $n$, and $k$ is illustrated by Figure \ref{fig:breakautos}. \medskip

\prf First, choose a $k\in\mathbb{N}$ that is larger than $1 + \frac{1}{\varepsilon}$. Then
\begin{equation}
\label{eq:k}
\frac{k-1}{k}>\frac{1}{1+\varepsilon}\,.
\end{equation}
Let $m\in\mathbb{N}$. By (\ref{eq:k}), there is an $n_0\in\mathbb{N}$ such that
\begin{equation}
\label{eq:n}
\forall\,n\geq n_0\colon(n-m)\cdot\frac{k-1}{k}\geq n\cdot\frac{1}{1+\varepsilon}+1\,.
\end{equation}
Let $n\in\mathbb{N}$ be sufficiently large, that is, $n\geq n_0$, and assume it fulfills (\ref{eq:1}). Then, the number of spheres $S_{v_0}(m+1), S_{v_0}(m+2),\ldots,S_{v_0}(n)$ that are not of the type $S_{v_0}(\lambda k+4)$, $\lambda\in\mathbb{N}$, is at least
\begin{equation}
\label{eq:spheres}
\left\lfloor(n-m)\cdot\frac{k-1}{k}\right\rfloor\geq\left\lfloor n\cdot\frac{1}{1+\varepsilon}+1\right\rfloor>\frac{n}{1+\varepsilon}\,\cdot
\end{equation}

Our goal is to $2$-color the vertices in these spheres in order to break all automorphisms that fix $v_0$ and act nontrivially on $B_{v_0}(m)$.

Let $\Aut(G,v_0)$ be the group of all automorphisms that fix $v_0$. Every $\phi\in\Aut(G,v_0)$ induces a permutation $\phi|B_{v_0}(n)$ of the vertices in $B_{v_0}(n)$. These permutations form a group $A$. If $\sigma$ and $\tau$ are different elements of $A$, then $\sigma\tau^{-1}\in A$ acts nontrivially on $B_{v_0}(n)$. By Lemma \ref{lem:main}, it also does so on $S_{v_0}(n)$, which means that $\sigma$ and $\tau$ do not agree on $S_{v_0}(n)$. Therefore, the cardinality of $A$ is at most $\left|S_{v_0}(n)\right|!$, for which the following rough estimate suffices for our purposes:
\begin{equation}
\label{eq:automorphisms}
\begin{array}{r@{\;}c@{\;}l} \left|S_{v_0}(n)\right|! & \leq & \left|S_{v_0}(n)\right|^{\left|S_{v_0}(n)\right|-1}\leq\left(\frac{n}{(1+\varepsilon)\log_2 n }\right)^{\frac{n}{(1+\varepsilon)\log_2 n }-1} \medskip \\ & \leq & n^{\frac{n}{(1+\varepsilon)\log_2 n }-1} = 2^{\left(\frac{n}{(1+\varepsilon)\log_2 n }-1\right)\log_2 n} \leq 2^{\frac{n}{1+\varepsilon}-1}\,. \end{array}
\end{equation}
It is clear that, if an element $\sigma\in A$ that acts nontrivially on $B_{v_0}(m)$ is broken by a suitable $2$-coloring of some spheres in $B_{v_{0}}(n)$, then all $\phi\in\Aut(G,v_0)$ with $\phi|B_{v_0}(n)=\sigma$ are broken at once. So it suffices to break all $\sigma\in A$ that act nontrivially on $B_{v_0}(m)$ by a suitable $2$-coloring of some spheres in $B_{v_{0}}(n)$ in order to ensure that all $\phi\in\Aut(G,v_0)$ that act nontrivially on $B_{v_0}(m)$ are broken.

Before doing this, let us remark that any element $\sigma\in A$ that acts nontrivially on the ball $B_{v_0}(m)$, also acts nontrivially on every sphere $S_{v_0}(m+1),\ldots,S_{v_0}(n)$. This is a consequence of Lemma \ref{lem:main}, and implies that we can break $\sigma$ by breaking the action of $\sigma$ on any one of the spheres $S_{v_0}(m+1),\ldots,S_{v_0}(n)$.

Now, consider the subset $S\subseteq A$ of all elements that act nontrivially on $B_{v_0}(m)$. As already remarked, every $\sigma\in S$ acts nontrivially on every sphere $S_{v_0}(m+1),\ldots,S_{v_0}(n)$. Hence, we can apply Lemma \ref{lem:halftwo} to break at least half of the elements of $S$ by a suitable coloring of $S_{v_0}(m+1)$. What remains unbroken is a subset $S'\subseteq S$ of cardinality at most $|S|/2$. Now, we proceed to the next sphere. We can break at least half of the elements of $S'$ by a suitable coloring of $S_{v_0}(m+2)$. What still remains unbroken, is a subset $S''\subseteq S$ of cardinality at most $|S|/4$.

Iterating the procedure, but avoiding spheres of the type $S_{v_0}(\lambda k+4)$, $\lambda\in\mathbb{N}$, we end up with the empty subset $\varnothing\subseteq S$ after at most $\log_{2}|S|+1\leq\log_{2}|A|+1\leq\frac{n}{1+\varepsilon}$ steps, see (\ref{eq:automorphisms}). This is less than the number of spheres not of the type $S_{v_0}(\lambda k+4)$, $\lambda\in\mathbb{N}$, between $S_{v_0}(m+1)$ and $S_{v_0}(n)$, see (\ref{eq:spheres}). Thus, we remain within the ball $B_{v_0}(n)$. Hence, all $s\in S$ and, therefore, all $\phi\in\operatorname{Aut}(G,v_{0})$ that act nontrivially on $B_{v_0}(m)$ are broken, and we are done. \qed

\begin{thm}
\label{thm:spheres}
Let $G$ be an infinite, locally finite, connected graph with infinite motion and $v_0\in\V(G)$. Moreover, let $\varepsilon>0$. If there exist infinitely many $n\in\mathbb{N}$ such that
\begin{equation}
\label{eq:smallspheres}
\left|S_{v_0}(n)\right| \leq \frac{n}{(1+\varepsilon)\log_2 n}\,,
\end{equation}
then the distinguishing number $\D(G)$ of $G$ is either $1$ or $2$.
\end{thm}

\prf Consider the $k\in\mathbb{N}$ provided by Lemma \ref{lem:breakautos}. First, we use Lemma \ref{lem:fixroot} to $2$-color all vertices in $B_{v_0}(k+3)$ and in $S_{v_0}(\lambda k+4)$, $\lambda\in\mathbb{N}$, such that, no matter how we color the remaining vertices, all automorphisms that move $v_{0}$ are broken.

Let $m_1=k+3$. Among all $n\in\mathbb{N}$ that satisfy (\ref{eq:smallspheres}) we choose a number $n_1\in\mathbb{N}$ that is larger than $m_1$ and sufficiently large to apply Lemma \ref{lem:breakautos}. Hence, we can $2$-color all vertices in $S_{v_0}(m_1+1),S_{v_0}(m_1+2),\ldots,S_{v_0}(n_1)$, except those in $S_{v_0}(\lambda k+4)$, $\lambda\in\mathbb{N}$, such that all automorphisms that fix $v_{0}$ and act nontrivially on $B_{v_0}(m_1)$ are broken. Next, let $m_2=n_1$ and choose an $n_2\in\mathbb{N}$ to apply Lemma \ref{lem:breakautos} again. Iteration of this procedure yields a $2$-coloring of $G$.

If an automorphism $\phi\in\Aut(G)\setminus\{\id\}$ moves $v_0$, then it is broken by our coloring. If it fixes $v_0$, consider a vertex $v$ with $\phi(v)\neq v$. Since $G$ is connected and $m_1<m_2<m_3<\ldots$, there is an $i\in\mathbb{N}$ such that $v$ is contained in $B_{v_0}(m_i)$. Hence, $\phi$ acts nontrivially on $B_{v_0}(m_i)$ and is again broken by our coloring. \qed

\begin{cor}
\label{cor:growth}
Let $G$ be an infinite, locally finite, connected graph with infinite motion and $v_0\in\V(G)$. Moreover, let $\varepsilon>0$. If there exist infinitely many $n\in\mathbb{N}$ such that
\begin{equation}
\label{eq:smallballs}
\left|B_{v_0}(n)\right| \leq \frac{n^2}{(2+\varepsilon)\log_2 n }\,,
\end{equation}
then the distinguishing number $\D(G)$ of $G$ is either $1$ or $2$. In particular, the Infinite Motion Conjecture holds for all graphs of growth $o(n^2/\log_2 n)$.
\end{cor}

\prf Let $n_1<n_2<n_3<\ldots$ be an infinite sequence of numbers that fulfill (\ref{eq:smallballs}). Notice that, for every $k\in\mathbb{N}$,
\begin{equation}
\sum_{i=1}^{n_k} \frac{i}{(1+\frac{\varepsilon}{2})\log_2 i } > \frac{{n_k}^2}{(2+\varepsilon)\log_2 n_k } \geq \left|B_{v_0}(n_k)\right| > \sum_{i=1}^{n_k} \left|S_{v_0}(i)\right|\,.
\end{equation}
Since
\begin{equation}
\lim_{k \rightarrow \infty}\left(\left(\sum_{i=1}^{n_k} \frac{i}{(1+\frac{\varepsilon}{2})\log_2 i }\right) - \frac{{n_k}^2}{(2+\varepsilon)\log_2 n_k }\right) = \infty\,,
\end{equation}
we infer that
\begin{equation}
\lim_{k \rightarrow \infty}\sum_{i=1}^{n_k}\left(\frac{i}{(1+\frac{\varepsilon}{2})\log_2 i }-\left|S_{v_0}(i)\right|\right) = \infty\,,
\end{equation}
and that, for infinitely many $i\in\mathbb{N}$,
\begin{equation}
\left|S_{v_0}(i)\right| < \frac{i}{(1+\frac{\varepsilon}{2})\log_2 i }\,.
\end{equation}
Hence, we can apply Theorem \ref{thm:spheres} to show that the distinguishing number $\D(G)$ of $G$ is either $1$ or $2$. \qedd

A result similar to Theorem \ref{thm:spheres} can also be obtained for graphs with countably many ends\footnote{Ends were first introduced by Freudenthal \cite{fr-45} in a topological setting, but here the definition of Halin \cite{ha-64} is more appropriate. For an accessible introduction to ends of infinite graphs see \cite{di-06}.}, none of which grows too fast. Readers not familiar with the notion of ends may safely skip the rest of this section, as the result is not used elsewhere in the paper.

\begin{thm}
Let $G$ be an infinite, locally finite, connected graph with countably many ends and infinite motion. Moreover, let $v_{0}\in\V(G)$ and $\varepsilon>0$. For an end $\omega$ of $G$ let $S_{v_0}^{\omega}(n)$ be the set of vertices in $S_{v_0}(n)$ that lie in the same connected component of $G \setminus B_{v_0}(n-1)$ as $\omega$. If, for every end $\omega$, there are infinitely many $n \in \mathbb N$ such that \begin{equation}
\label{eq:smallomegaspheres}
\left|S_{v_0}^\omega(n)\right| \leq \frac{n}{(1+\varepsilon)\log_2 n}\,,
\end{equation}
then the distinguishing number $\D(G)$ of $G$ is either $1$ or $2$.
\end{thm}

\prf Basically the proof consists of three steps. First we color part of the vertex set in order to break all automorphisms that move $v_0$. In the second step we break all automorphisms in $\Aut (G,v_0)$ that do not fix all ends of the graph by coloring some other vertices. Finally, we color the remaining vertices to break the rest of the automorphisms.

In order to break all automorphisms that move $v_0$ we apply Lemma~\ref{lem:fixroot}, just as in the proof of Theorem~\ref{thm:spheres}. The only difference is that we choose $k$ twice as large as proposed by Lemma~\ref{lem:breakautos}, because we would like to color some additional spheres in the second step of the proof before applying an argument similar to that in Lemma~\ref{lem:breakautos}.

For the second step consider the spheres $S_{v_0}(\frac{2\lambda+1}{2}k+4)$, $\lambda\in\mathbb{N}$. We wish to color those spheres such that every automorphism that fixes $v_0$ and preserves the coloring also fixes every end of $G$.

It is not hard to see that the sets $S_{v_0}^\omega(\frac{2\lambda+1}{2}k+4)$, $\omega$ an end of $G$, $\lambda\in\mathbb{N}$, carry the following tree structure. Consider $v_0$, the root, which is connected by an edge to $S_{v_0}^\omega(\frac{3}{2}k+4)$ for each end $\omega$. For every end $\omega$ of $G$ and every $\lambda\in\mathbb{N}$, draw an edge from $S_{v_0}^\omega(\frac{2\lambda+1}{2}k+4)$ to $S_{v_0}^\omega(\frac{2\lambda+3}{2}k+4)$. To see that this is indeed a tree just notice that if $S_{v_0}^{\omega_1}(n)=S_{v_0}^{\omega_2}(n)$, then, for every $m \leq n$, $S_{v_0}^{\omega_1}(m)=S_{v_0}^{\omega_2}(m)$. So there cannot be any circles. By construction, this tree structure is infinite, locally finite, and does not have any endpoints.

Next, notice that every automorphism $\phi \in \Aut (G,v_0)$ that does not fix all ends also acts as an automorphism on this tree structure. By \cite{wazh-07}, the distinguishing number of infinite, locally finite trees without endpoints is at most $2$. Therefore it is possible to $2$-color the sets $S_{v_0}^\omega(\frac{2\lambda+1}{2}k+4)$, $\omega$ an end of $G$, $\lambda\in\mathbb{N}$, such that every such automorphism is broken. It is also worth noting that so far we did not use the countability of the end space of $G$, nor did we use the growth condition on the ends.

Let us turn to the third step of the proof. So far we have colored the ball $B_{v_{0}}(k+3)$ and the spheres $S_{v_0}(\frac{\lambda}{2}k+4)$, $\lambda\geq 2$, in a way that color preserving automorphisms fix $v_0$ and move every $S_{v_0}^\omega (n)$ into itself. Consider such an automorphism $\phi$, which acts nontrivially on $G$. If we remove the fixed points of $\phi$ from $G$, then the infinite motion of $G$ implies that the resulting graph has only infinite components. Hence, there is a ray in $G$ which contains no fixed point of $\phi$. The image of this ray must lie in the same end $\omega$. Thus, there is an index $n_0$, such that, for every $n \geq n_0$, $\phi$ acts nontrivially on $S_{v_0}^{\omega}(n)$.

Let $(\omega_i)_{i \in {\mathbb N}}$ be an enumeration of the ends of $G$. Choose a function $f\colon {\mathbb N} \to {\mathbb N}$ such that, for every $i\in\mathbb{N}$, $f^{-1}(i)$ is infinite. Assume that all spheres up to $S_{v_0}(m)$ have been colored in the first $i-1$ steps. In the $i$-th step we would like to color some more spheres in order to continue breaking all automorphisms in $\Aut(G, v_0)$ that act nontrivially on each of the spheres $S_{v_0}^{\omega_{f(i)}}(n)$, $n\geq m$. This can be done by exactly the same argument as the one used in the proof of Lemma~\ref{lem:breakautos}.

As we already mentioned, every automorphism that was not broken in the first two steps acts by nontrivially on the rays of some end. Since, in the procedure described above, every end is considered infinitely often, it is clear that every such automorphism will eventually be broken. This completes the proof. \qed

\section{Graphs with higher cardinality}
\label{sec:generalization}
If a graph $G$ has trivial automorphism group, then $G$ is obviously $1$-distinguishable, that is, $\D(G) = 1$. From now on we assume that our graphs $G$ have nontrivial automorphism group. In this case, the \emph{motion} $\m(G)$ of $G$ is defined as
\begin{equation}
\m(G)=\min_{\phi\in\Aut(G)\setminus\{\id\}}\m(\phi)\,.
\end{equation}
As already mentioned, the Infinite Motion Conjecture does not hold for graphs of higher cardinality. An example is the Cartesian product $G = K_{\textswab n} \cp K_{\textswab m}$ of two complete graphs on infinitely many vertices $\textswab{n}$ and $\textswab{m}$ with $2^\textswab n < \textswab m$. By \cite{imkltr-07}, $G$ has motion $\textswab n$, but $\D(G) > \textswab n$.

The question arises whether one can adapt the Infinite Motion Conjecture to graphs of higher cardinality. The starting point is \cite[Theorem 1]{rusu-98}. It asserts that a finite graph $G$ is $2$-distinguishable if $\m(G) > 2\log_2|\!\Aut(G)|$. However, a second look at the proof shows that the inequality sign can be replaced by $\geq$. For details see Section \ref{se:rusu}. For finite graphs we thus infer that
\begin{equation}
\label{eq:russell}
\m(G) \geq 2\log_2|\!\Aut(G)|\quad\text{implies}\quad\D(G)=2,
\end{equation}
which can also be written in the form
\begin{equation*}
|\!\Aut(G)|\leq 2^{\frac{\m(G)}{2}} \quad\text{implies}\quad \D(G)=2\,.
\end{equation*}
Notice that $2^{\frac{\m(G)}{2}} = 2^{\m(G)}$ if $\m(G)$ is infinite. We are thus tempted to conjecture for graphs $G$ with infinite motion that $ |\!\Aut(G)| \leq 2^{\m(G)}$ implies $ D(G)=2$. We formulate this conjecture as the
\bigskip

{\bf Motion Conjecture.} \emph{Let $G$ be a connected graph with infinite motion $\m(G)$ and} $|\!\Aut(G)| \leq 2^{\m(G)}\,.$ \emph{Then $\D(G)=2$.}
\bigskip

How does this compare with the Infinite Motion Conjecture? It asserts that the distinguishing number of a locally finite, connected graph $G$ is $2$ if $\m(G)$ is infinite. Since locally finite graphs are countable, the condition that $\m(G)$ is infinite is equivalent to $\m(G) = \aleph_0$. Furthermore, for countable graphs we have
\[ |\!\Aut(G)| \leq \aleph_0^{\aleph_0} = 2^{\aleph_0}.\]% = \textswab{c}\,. \]
Hence, for countable graphs, and thus also for locally finite, connected graphs with infinite motion, the inequality of the Motion Conjecture is automatically satisfied, which means that the Infinite Motion Conjecture is a special case of the Motion Conjecture.

Now, let us focus on the two results from \cite{imsmtuwa-xx} that imply the validity of the Infinite Motion Conjecture for graphs with countable group.

\begin{thm}
\label{thm:simon}
Let $G$ be a locally finite, connected graph that satisfies $\aleph_0 \leq |\!\Aut(G)| < 2^{\aleph_0}$. Then $|\!\Aut(G)| = {\aleph_0}$, $\m(G) = \aleph_0$, and $\D(G) = 2$.
\end{thm}

Notice that the only thing that is required here, besides local finiteness and connectedness, is an upper and a lower bound on the size of $\Aut(G)$. And it turns out, that $\Aut(G)$ is countable, even without the continuum hypothesis. Even infinite motion and $\D(G) = 2$ are consequences of this restriction on the size of the automorphism group.

\begin{thm}
\label{thm:tucker}
Let $G$ be a countably infinite, connected graph that satisfies the conditions $|\!\Aut(G)| \leq \m(G)$ and $\m(G) = \aleph_0$. Then $\D(G)=2$.
\end{thm}

Here, without local finiteness, one cannot drop the assumption of infinite motion. If we assume that $\Aut(G)$ has smaller cardinality than the continuum, then we can ensure $2$-distinguishability if the continuum hypothesis holds, but we do not know whether this is really necessary.

\begin{cor}
\label{cor:tucker}
Let $G$ be a countably infinite, connected graph with infinite motion. If the continuum hypothesis holds, and if $ |\!\Aut(G)| < 2^{\m(G)}$, then $D(G) =2$.
\end {cor}

The next theorem shows that Theorem \ref{thm:tucker} also holds for graphs of higher cardinality and uncountable motion.

\begin{thm}
\label{thm:special}
Let $G$ be a connected graph with uncountable motion. Then $|\!\Aut(G)| \leq \m(G)$ implies $\D(G)=2$.
\end{thm}

\prf Set $\textswab{n} = |\!\Aut(G)|$, and let $\zeta$ be the smallest ordinal number whose underlying set has cardinality $\textswab{n}$. Furthermore, choose a well ordering $\prec$ of $A = \Aut(G)\setminus \{\id\}$ of order type $\zeta$, and let $\alpha_0$ be the smallest element with respect to $\prec$. Then the cardinality of the set of all elements of $A$ between $\alpha_0$ and any other $\alpha \in A$ is smaller than $\textswab{n} \leq \m(G)$.

Now we color all vertices of $G$ white and use transfinite induction to break all automorphisms by coloring selected vertices black.

{\sc Induction base\,} By the assumptions of the theorem, there exists a vertex $v_0$ that is not fixed by $\alpha_0$. We color it black. This coloring breaks $\alpha_0$.

{\sc Induction step\,} Let $\beta \in A$. Suppose we have already broken all $\alpha \prec \beta$ by pairs of distinct vertices $(v_\alpha, \alpha(v_\alpha))$, where $v_\alpha$ is black and $\alpha(v_\alpha)$ white. Clearly, the cardinality of the set $R$ of all $(v_\alpha, \alpha(v_\alpha))$, $\alpha \prec \beta$, is less than $\m(G) \geq \textswab{n}$. By assumption, $\beta$ moves at least $\m(G)$ vertices. Since there are still $\textswab{n}$ vertices not in $R$, there must be a pair of vertices $(v_{\beta}, \beta(v_\beta))$ that does not meet $R$. We color $v_{\beta}$ black. This coloring breaks $\beta$. \qed

\begin{cor}
\label{cor:special}
Let $G$ be a connected graph with uncountable motion. If the general continuum hypothesis holds, and if $ |\!\Aut(G)| < 2^{\m(G)}\!,$\, then $D(G) =2$.
\end {cor}

\prf Under the assumption of the general continuum hypothesis $2^{\m(G)}$ is the successor of $\m(G)$. Hence $|\!\Aut(G)| \leq \m(G) $, and the assertion of the corollary follows from Theorem \ref{thm:special}. \qed

\section{The Motion Lemma of Russell and Sundaram}
\label{se:rusu}

In order to show that a finite graph $G$ is $2$-distinguishable if $\m(G) > 2\log_2|\!\Aut(G)|$, Russell and Sundaram \cite{rusu-98} first defined the cycle norm of an automorphism $\phi$. If
\[\phi = (v_{11}v_{12} \ldots v_{1l_1})(v_{21} \ldots v_{2l_2}) \ldots (v_{k1}\ldots v_{kl_k} )\,,\]
then the \emph{cycle norm} $\cn(\phi)$ of $\phi$ is
\[ \cn(\phi) = \sum_{i=1}^k(l_i -1)\,. \]
The cycle norm $\cn(\phi)$ is related to graph distinguishability as follows: Let $G$ be randomly $2$-colored by independently assigning each vertex a color uniformly from $\{\text{black},\text{white}\}$. Then the probability that every cycle of $\phi$ is monochromatic is $2^{-\cn(\phi)}$. In this case, $\phi$ preserves the coloring so chosen.

Further, they define the cycle norm $\cn(G)$ of a graph $G$ as
\[ \cn(G) = \min_{\phi \in \Aut(G)\setminus\{\id\}} \cn(\phi)\,. \]
We now reprove Theorem 2 of \cite{rusu-98} with $\geq$ instead of $>$. Because $\cn(\phi)\geq\m(\phi)/2$ and thus $\cn(G)\geq\m(G)/2$ we infer from Theorem \ref{thm:rusunew} below that $G$ is $2$-distinguishable if $\m(G)\geq 2\log_2|\!\Aut(G)|$. We propose to call this result ``Motion Lemma of Russell and Sundaram''. Actually, the only difference from the original proof is the insertion of the middle term in (\ref{eq:sund}).

\begin{thm}
\label{thm:rusunew}
Let $G$ be a finite graph, and $\cn(G)\log d \geq \log |\!\Aut(G)|$. Then $G$ is $d$-distinguishable, that is, $\D(G)\leq d$.
\end{thm}

\prf Let $\chi$ be a random $d$-coloring of $G$, the probability distribution being given by selecting the color of each vertex independently and uniformly in the set $\{1, \ldots, d\}$. For a fixed automorphism $\phi\in\Aut(G)\setminus\{\id\}$ consider the probability that the random coloring $\chi$ is preserved by $\phi$:
\begin{equation}
\Prob_\chi[\forall\,v\,:\,\chi(\phi(v)) = \chi(v)] = \left(\frac{1}{d}\right)^{\cn(\phi)} \leq \left(\frac{1}{d}\right)^{\cn(G)}\,.
\end{equation}
Collecting these events yields the inequality
\begin{equation}
\label{eq:sund}
\begin{array}{r@{\;}c@{\;}l} \Prob_\chi[\exists\,\phi\in\Aut(G)\setminus\{\id\}\,\forall\,v\,:\,\chi(\phi(v)) = \chi(v)] & \leq & (|\!\Aut(G)|-1)\left(\frac{1}{d}\right)^{\cn(G)} \medskip \\ & < & |\!\Aut(G)|\left(\frac{1}{d}\right)^{\cn(G)}\,. \end{array}
\end{equation}
By hypothesis the last term is at most 1. Thus there exists a coloring $\chi$ such that, for every $\phi\in\Aut(G)\setminus\{\id\}$, there is a $v$ for which $\chi(\phi(v))\neq\chi(v)$, as desired. \qedd

{\bf Acknowledgement.} We thank the referees for their comments and remarks, as they contributed considerably to the readability of the paper. Furthermore, we are grateful to Norbert Sauer and Claude Laflamme for their suggestions pertaining to Theorem \ref{thm:special} and Corollary \ref{cor:special}.

\end{document}